\documentclass[journal]{IEEEtran}

\ifCLASSINFOpdf
\else
   \usepackage[dvips]{graphicx}
\fi
\usepackage{url}

\usepackage{graphicx}
\usepackage{hashem}

\begin{document}

\title{Supervised learning of analysis-sparsity priors with automatic differentiation}

\author{Hashem Ghanem, Joseph Salmon, Nicolas Keriven, and Samuel Vaiter
  \thanks{The authors acknowledge the support of ANR Grava ANR-18-CE40-0005, ANR GRandMa ANR-21-CE23-0006 and ANER BFC RAGA.}
	\thanks{HG is with CNRS, IMB, Univ. de Bourgogne; JS is with Univ. Montpellier; NK is with CNRS, GIPSA-lab; SV is with CNRS, LJAD, Univ. Côte d'Azur.}\vspace{-1cm}}

\markboth{}
{}
\maketitle

\begin{abstract}
	Sparsity priors are commonly used in denoising and image reconstruction. For analysis-type priors, a dictionary defines a representation of signals that is likely to be sparse.
	In most situations, this dictionary is not known, and is to be recovered from pairs of ground-truth signals and measurements, by minimizing the  reconstruction error.
	This defines a hierarchical optimization problem, which can be cast as a bi-level optimization.
	Yet, this problem is unsolvable, as reconstructions and their derivative \wrt the dictionary have no closed-form expression.
	However, reconstructions can be iteratively computed using the Forward-Backward splitting (FB) algorithm.
	In this paper, we approximate reconstructions by the output of the aforementioned FB algorithm.
	Then, we leverage automatic differentiation to evaluate the gradient of this output \wrt the dictionary, which we learn with projected gradient descent.
	Experiments show that our algorithm successfully learns the 1D Total Variation (TV) dictionary from piecewise constant signals.
	For the same case study, we propose to constrain our search to dictionaries of 0-centered columns, which removes undesired local minima and improves numerical stability.
\end{abstract}

\begin{IEEEkeywords}
	Sparsity, dictionary learning, total variation, bi-level optimization, automatic differentiation.
\end{IEEEkeywords}

\IEEEpeerreviewmaketitle

\section{Introduction}
\IEEEPARstart{D}{enoising}  is a widely-tackled problem that emerges in many fields,
where the goal is to restore signals from  noisy observations. This includes that the model of the imaging system is \emph{a priori} known: $\V{y= w+\varepsilon}$, where $\V{w}, \V{y} \in \bbR^p$ are the true and the measured signals, respectively, and $\V{\varepsilon}\in \bbR^p$ is additive noise.
In addition, a prior hypothesis on the nature of the signal, or an embedding of it, might be available, like sparsity \cite{mccann2020supervised}. Such extra knowledge can be incorporated in the optimization process to get better reconstructions (\eg a higher Signal to Noise Ratio (SNR)).

Sparsity priors exist in two forms \cite{elad2007analysis}:
\emph{i)} synthesis (traditional) sparsity where $\V{w}= \V{D}\V{u}$, $\V{D}$ is a linear operator, and $\V{u}$ is sparse;
\emph{ii)} analysis sparsity where $\V{v} = \V{D}^\top \V{w} \in \bbR^{m}$ is sparse. This paper is interested in the latter.
As a convex surrogate, this prior is enforced on $\V{w}$ by adding the term $\|\V{D^\top w}\|_1$ to the (generally quadratic) loss function \cite{Mancera2006L0}, where $\| \cdot \|_1 $ is the $\ell_1$ norm. Putting all together,  the problem consists in finding:
$\V{\hat{w}} = \argmin_{\V{w}} \|\V{y} - \V{w}\|_2^2 + \lambda\| \V{D}^\top \V{w}\|_1$
for some regularization amplitude parameter $\lambda \geq 0$.
The linear operator $\V{D}$ is either user-defined, or learnt directly from data.

The main problem we tackle in this work is to extract both $\V{D}$ and $\lambda$ from data with supervised learning.
Having that $\lambda\| \V{D}^\top \V{w}\|_1 = \| (\lambda\V{D})^\top \V{w}\|_1$, this problem is equivalent to extracting the product $\lambda\V{D}$ as one object, thus we stop notating $\lambda$  explicitly from now on and keep $\V{D}$.
To avoid trivial solutions and undesired local minima, it is common to restrict the search space to an admissible set $\calC \subset \bbR^{p\times m}$, where dictionaries have a specific property. In \cite{ravishankar2015sparsifying} for instance, $\V{D}$ is forced to be orthogonal, \ie $\calC = \{\V{D}; \V{DD}^\top= \V{I}_p\}$. In \cite{peyre2011learning}, $ \V{D}$ is constrained to be a convolution dictionary.
However, it is still challenging to find an admissible set that performs well in all applications of this problem \cite{yaghoobi2012noise}.
In this work, we don't commit to find such universal set.
Formally, we suppose we have a dataset  $(\V{y}_l, \V{w}_l)_{l=1}^L$, that includes $L$ pairs of measurements and their associated ground truth signals.
Knowing such pairs, our task is to find the operator $\V{D}$ that minimizes the mean squared error between reconstructions and ground truth signals:
\begin{subequations}\label{eq:bilevel_problem}
	\begin{equation}
		\label{eq:D_learning_opt_cost}
		\hat{\V{D}}\in
		\argmin_{\V{D}\in \bbR^{p\times m}} \calE (\V{D})
		= \sum_{l=1}^L \big \Vert \hat{\V{w}}\left(\V{D}, \V{y}_l\right) - \V{w}_l \big\Vert_2^2 + \iota_\calC(\V{D})
	\end{equation}
	\st
	\begin{equation}
		\label{eq:denoising_opt_cost}
		\hat{\V{w}}(\V{D}, \V{y}) =
		\argmin_{\V{w}\in \bbR^p} \frac{1}{2} \| \V{y} - \V{w} \|_2^2 + \| \V{D}^\top \V{w}\|_1 \enspace,
	\end{equation}
\end{subequations}
where $\calC$ is an admissible set of dictionaries, and $\iota_\calC$ is the indicator function: $\iota_\calC(\V{D})=0$ if $\V{D}\in \calC$ and $+\infty$ otherwise.

One can recast \Cref{eq:bilevel_problem} as a \emph{bilevel optimization problem}: in the outer problem, we learn the model $\V{D}$ as in \cref{eq:D_learning_opt_cost}  while in the inner problem, we denoise measurements following \cref{eq:denoising_opt_cost}.
We already know that the inner part can be solved applying the Forward-Backward splitting (FB) algorithm on the dual problem \cite{chambolle2010introduction}.
However, due to the $\ell_1$ norm, neither the solution nor its gradient \wrt $\V{D}$ have closed-form expressions.
Thus, the minimizer $\hat{\V{D}}$ cannot be derived analytically nor obtained  with gradient-based methods.

To illustrate how our algorithm can recover a dictionary $\hat{\V{D}}$, we consider the well-known problem of \emph{1D piecewise constant signals reconstruction} as a case-study \cite{chambolle2010introduction}, where this prior indicates that $(w_2-w_1,\dots,w_1-w_p)^\top$ is sparse.
The estimator is often written as an instance of \cref{eq:denoising_opt_cost}, with $\V{D}=\V{D}_{TV}$ is the dictionary associated to the 1D Total Variation (TV) regularization: for all $i\in \{1,\dots,p\}; \V{D}_{i,i}=-1, \V{D}_{i+1,i}=1$, and 0 otherwise, up to rescaling.

\textbf{Contribution}\quad We approximately recover the analysis-sparsity operator $\hat{\V{D}}$ by:
\emph{i)} replacing the true minimizer $\hat{\V{w}}(\V{D}, \V{y})$ by the output of the FB algorithm applied on the dual problem;
\emph{ii)} deploying \emph{automatic differentiation}, a technique  capable of evaluating the gradient of an algorithm \wrt input variables, to solve \cref{eq:D_learning_opt_cost} with projected gradient descent.
We empirically show that our method recovers the TV dictionary $\V{D}_{TV}$ from piecewise constant signals.
Finally, for the same case study, we reduce the admissible set $\calC$ to dictionaries with columns summing up to zero, and empirically prove that this increases stability and extracts the TV operator with higher quality than previous methods.

\textbf{Related work} \quad The bilevel problem was first posed in \cite{peyre2011learning}, where the authors smoothed the $\ell_1$ norm so that the derivative of $\hat{\V{w}}(\V{D},\V{y})$ \wrt $\V{D}$ has a closed-form expression.
Then, they applied gradient descent to find a local minimizer $\V{D}$.
However, the sought for sparsity is degraded with this methodology.
In \cite{sprechmann2013supervised}, a relaxation regime of $\ell_1$ with the smooth $\ell_2$ norm is adopted, with a relaxation parameter easy to assign.
Recently in \cite{mccann2020supervised}, a formula of $\hat{\V{w}}(\V{D},\V{y})$ is obtained under some conditions. This formula includes inverting a large matrix, which is computed iteratively while Automatic Differentiation (AD) is used to evaluate gradients.
In \cite{yaghoobi2011analysis}, and without a proof of convergence, they  restrict $\V{D}$ to have unit columns norm, that also satisfy $\V{D}\V{D}^\top=\V{I}$, in order to avoid trivial solutions (\eg repeated columns).
The work in \cite{chambolle2020learning} solves the 2D piecewise  signals reconstruction problem using AD, constrained by learning \emph{convolution-type} dictionaries of kernels with small support, which improves the quality of standard 2D TV. Such strong constraints are not considered in our work, however we will see that simple column-centering suffices to learn a high-quality dictionary in the 1D version of this problem.

\section{Automatic Differentiation (AD)}
\label{sec:preliminaries}

To compute gradients, one often manually writes down its analytical expression.
Yet, this can only be performed for functions with closed-form expression, which is not the case in \cref{eq:denoising_opt_cost}.
An alternative is symbolic differentiation, automatically performed by computer tools like Mathematica \cite{grabmeier2003computer} when dealing with syntax tree expressions.
A third approach consists in approximating gradients using numerical differentiation. It is easy to implement though exposed to round-off errors \cite{jerrell1997automatic} and expensive in case of a high number of variables (here $p \times m$).

Automatic Differentiation (AD), that is of interest in this work, mitigates the previous drawbacks.
AD manipulates computation flow in a computer program, with all numerical computations reduced to compositions of elementary operations, for which derivative rules are known \cite{verma2000introduction}.
In such computer programs, a value must be assigned to each input variable, thus, any operation in the code will result in a variable with a numerical value.
During the execution of the program, AD consists in \cite{baydin2018automatic}: 1) keeping a trace to all intermediate variables that are dependent on the ones we want to differentiate for; 2) once an operation is performed on a dependent variable, say $v_i$, to evaluate another $v_j$, directly computing the derivative value $dv_j/dv_i$; 3) accumulating derivatives in step 2 through the chain rule.
This gives the derivative value of the whole composition \wrt a chosen variable.
AD can trace numerical computations in recursion functions and in control-flow (\emph{if, while, for}) statements. Thus, AD efficiently differentiates not only closed-form formulas, but also implemented algorithms.
This makes AD suitable for gradient-based optimization.
AD can be implemented in two ways: \emph{forward} mode and \emph{reverse} mode.
In our work, we adopt the reverse mode\footnote{this is the mode implemented in the PyTorch framework}.
Let us assume having $y= (y_1,\dots, y_m)$ as a function of the variable $x=(x_1,\dots,x_n)$, with $v$ as an intermediate variable.
The backward scheme takes in input the vector $\bar{\V{y}} = (\bar{y}_1, \dots, \bar{y}_m)$, and accumulates gradients as follows: first $\bar{\V{v}}= \V{J}_y^\top(v) \bar{\V{y}}$, then   $\bar{\V{x}} =\V{J}_v^\top(x) \bar{\V{v}}$, where $\left(\V{J}_y(v)\right)_{i,j}= dy_i/dv_j$. The output $\bar{\V{x}}$ is nothing but the transpose Jacobian-vector product $\V{J}^\top_y(x) \bar{\V{y}}$.

\begin{figure*}
	\centering
	\includegraphics[width=\textwidth]{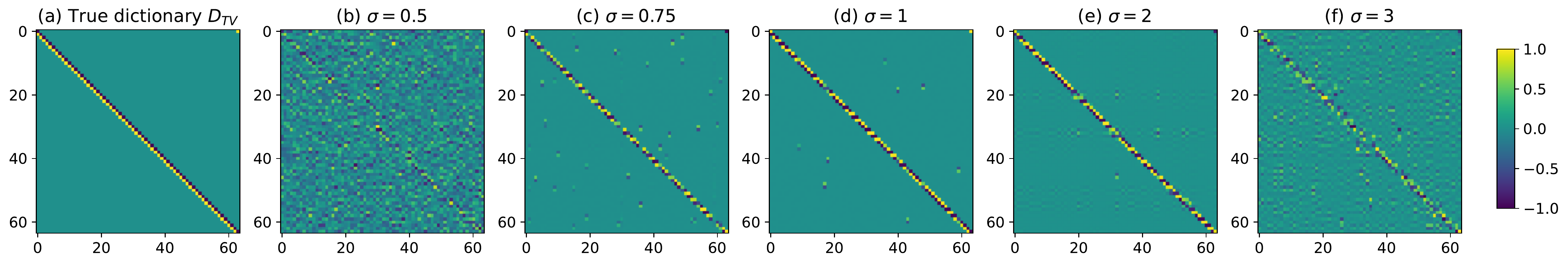}
	\caption{Performance of our projected gradient descent algorithm \ref{alg},  \wrt to noise level. We plot a sorted view of the dictionary $\hat{\V{D}}$.  We rescale all dictionaries to $[-1,1]$ for a better visualization of the structure recovered in $\hat{\V{D}}$ . (a) $\V{D}_{TV}$ our algorithm is expected to learn. From (b) to (f): $\hat{\V{D}}$ for different values of $\sigma$, the standard deviation of noise. Jumps in the dataset have the same amplitude: $0 \rightarrow  10$  or  $10 \rightarrow  0$.}
	\label{fig:varying_sigma}
\end{figure*}

\section{Proposed algorithm}
\label{sec:proposed_algorithm}

We solve the \emph{dual} problem of \cref{eq:denoising_opt_cost} with Forward-Backward splitting (FB).
Using  Automatic Differentiation, we obtain gradients of the previous FB algorithm \wrt $\V{D}$.
These gradients are used to learn a local minimum $\hat{\V{D}}$ using gradient descent, while projecting $\V{D}$ on the admissible set $\calC$ at every iteration.

\subsection{Deriving the dual problem of \cref{eq:denoising_opt_cost}}
\label{sub:Deriving_the_ dual_problem}

The term $\| \V{D}^\top \V{w}\|_1$ in \cref{eq:denoising_opt_cost} is neither differentiable \wrt $\V{w}$, nor has a simple proximal operator that can be efficiently performed.
Hence, we cannot apply the gradient descent or the FB algorithm directly to evaluate $\hat{\V{w}}(\V{D}, \V{y})$.
However, the latter is possible if we tackle this optimization from the dual perspective \cite{rockafellar1974conjugate}.
In fact, one can prove that \cref{eq:denoising_opt_cost} is equivalent to its dual problem \cite{chambolle2010introduction}:

\begin{equation}
	\label{eq:dual_problem}
	\hat{\V{z}}(\V{D}, \V{y}) =
	\argmin_{\V{z}\in \bbR^m} ~ \frac{1}{2} \big\Vert \V{Dz}- \V{y} \big\Vert_2^2
	~+~ \iota_{B_\infty}(\V{z})\enspace,
\end{equation}
where $B_\infty = \{\V{z}\in \bbR^m ; ~| z_i | \leq 1, \ \forall i\in[m]\}$ is the unit ball of the $\ell_\infty$ norm.
The target recovery $\hat{\V{w}}$ is then given by:
\begin{equation}
	\label{eq:dual_to_primal}
	\hat{\V{w}}(\V{D}, \V{y})
	=  \V{y} -\V{D}~\hat{\V{z}}(\V{D}, \V{y}).
\end{equation}

\subsection{Solving the dual problem with the  FB algorithm}
\label{sub:Solving_the_dual_problem_with_FB}

Fortunately, the term $\iota_{B_\infty}(\V{z})$ has a simple proximal function given by $\Pi_{B_\infty}$: the orthogonal projection on the ball $B_\infty$. So indeed, we solve \cref{eq:dual_problem} with an accelerated FB algorithm, namely the Fast Iterative Shrinkage-Thresholding Algorithm (FISTA) \citep{beck2009fast}.
At each iteration we update $\V{z}$ as follows:
\begin{align}
	\label{eq:dual_FB_update}
	\begin{split}
		\V{z}_{i+1} &=
		\prox_{\iota_{B_\infty}} \big(\V{q}_i -\eta_1 \nabla (\frac{1}{2} \big\Vert \V{D}\V{q}_i- \V{y} \big\Vert_2^2) \big)  \\
		&=\Pi_{B_\infty} \big(\V{q}_i - \eta_1 \V{D}^\top (\V{D}\V{q}_i - \V{y}) \big)
	\end{split}
\end{align}
\st $\V{q}_1 = \V{z}_0$ is the initialization of $\V{z}$, $\eta_1$ is the step size, and:
\begin{align}
	\label{eq:update_q}
	\begin{split}
		\V{q}_{i+1} &=
		\V{z}_i + \frac{t_i-1}{t_{i+1}} (\V{z}_i - \V{z}_{i-1})\\
		t_{i+1} &=
		\frac{1}{2} (1+ \sqrt{1+4t_i^2 })	\enspace; t_1 = 1.
	\end{split}
\end{align}

To conclude with the inner part: having the matrix $\V{D}$ and the vector $\V{y}$ as \emph{input}, one can compute $\hat{\V{z}}(\V{D}, \V{y})$ with a sufficient number of updates as in \cref{eq:dual_FB_update}, then get $\hat{\V{w}}(\V{D}, \V{y})$ \emph{in output} using \cref{eq:dual_to_primal}.

\subsection{Outer loop design to learn the dictionary}
\label{sub:Outer_loop_design}

Let us first notice that when the admissible set $\calC=\bbR^{p\times m}$, the solution $\hat{\V{D}}$ is \emph{not unique} as $\Vert \cdot \Vert_1$  is invariant to the coefficients order in a vector, \eg $\Vert (u_1, u_2)^\top \Vert_1 = \Vert (u_2, u_1)^\top \Vert_1$. Thus, permuting columns in $\V{D}$ will lead to the same $\hat{\V{w}}(\V{D}, \V{y})$ in \cref{eq:denoising_opt_cost}, which means the same cost $\calE(\V{D})$. In this paper, we consider any solution that is a minimizer of $\calE(\V{D})$. We  also assume that the second dimension of $\V{D}$ is given, while optimizing for it might be the subject of a future work.

Towards our goal, we use AD to get the gradient of the mean squared error (MSE) term in $\calE(\V{D})$, by tracing the FB algorithm that solves \cref{eq:denoising_opt_cost}. In fact, to reduce computational cost, we do not compute the full MSE but sample a batch of training signals at each iteration (see Alg.~\ref{alg}). We denote this term by $\calE_{\rm MSE}(\V{D})$. Since proving the convergence of AD's Jacobian to the variational one is complicated in such non-smooth setting, we replace the true reconstruction $\hat{\V{w}}(\V{D},\V{y})$ with the output produced by the FB algorithm, and empirically show that it is a good proxy. For simplicity, we keep the same notation for this output.
We randomly initialize $\V{D}_0$, then we start each iteration $t$ by setting \emph{PyTorch AD framework}  to record operations on $\V{D}_t$.
We can then compute the output $\hat{\V{w}}(\V{D}_t, \V{y})$, and $\calE_{\rm MSE}(\V{D}_t)$.
Now, we use the PyTorch AD to get the gradient $\nabla \calE_{\rm MSE}(\V{D}_t)$, and we update the estimated dictionary as
\begin{equation}
	\label{eq:outer_loop_gradient_step}
	{\V{D}}_{t+1} =
	\V{D}_t- \eta_2 \nabla \calE_{\rm MSE}(\V{D}_t)\enspace,
\end{equation}
where $\eta_2>0$ is a step size. Lastly, we project $\V{D}_t$ on the admissible set $\calC$ by computing $\Pi_\calC(\V{D}_t)$.
We obtain Algorithm~\ref{alg}.

\begin{algorithm}[h]
	\DontPrintSemicolon
	\KwInput{$\left\{ (\V{w}_l,\V{y}_l) \right\}_{l\in\{1, \dots , L\}}$: dataset.}
	\KwOutput{$\hat{\V{D}}$: minimizer to $\calE (\V{D})$ in \cref{eq:D_learning_opt_cost}.}
	\Hyp{$m, \eta_1,\eta_2, max\_itr1, max\_itr2, batch\_sz$\\}
	\Algo{\\}
	Initialize $\V{D}$  \emph{iid} from $\calN(0,10^{-4})$.\\
	Set PyTorch AD to track computations on $\V{D}$.\\
	\For{$t$ in $\{0,\dots,max\_itr2-1\}$:\\}{
		$\calE_{MSE}(\V{D})\leftarrow 0$\\
		\For{$l \in \{t*batch\_sz,\dots,(t+1)*batch\_sz-1\}$\\}{
			Initialize $\V{q}$ \emph{iid} from $\calN(0,1)$.\\
			\For{$i$ in $\{1,\dots,max\_itr1\}$:\\}{
				$\V{z}\leftarrow \Pi_{B_\infty} \big(\V{q} - \eta_1 \V{D}^\top (\V{D}\V{q} - \V{y}_l) \big)$.\\
				Update $\V{q}$ as in \cref{eq:update_q}.\\
			}
			$\hat{\V{w}}(\V{D}, \V{y}_l) \leftarrow  \V{y}_l -\V{D}\V{z}.$\\
			$\calE_{MSE} (\V{D}) \leftarrow \calE_{MSE} (\V{D}) + \big \Vert \hat{\V{w}}\left(\V{D}, \V{y}_l\right) - \V{w}_l \big\Vert_2^2.$\\
		}
		$\nabla \calE_{MSE} (\V{D})\leftarrow$ PyTorch AD gradient.\\
		$	{\V{D}} \leftarrow \Pi_\calC\big(\V{D}- \frac{ \eta_2}{batch\_sz}\nabla \calE_{MSE} (\V{D})\big)$.\\
	}
	Return $\hat{\V{D}}=\V{D}$

	\caption{AD-based projected gradient descent}\label{alg}
\end{algorithm}

\section{Experiments}
\label{sec:Experiments}

We consider the 1D piecewise constant signals reconstruction problem, with $p=m=64$ in all experiments. Ground-truth examples $\V{w}$ are generated \st they have 4 discontinuities, thus 4 constant pieces.
The coefficients where discontinuities take place are randomly chosen in each $\V{w}$. Their amplitude is either fixed ($0 \rightarrow  10$  or  $10 \rightarrow  0$), or randomly sampled \st they happen between two values in $[0,10]$, to be mentioned when necessary.
Observations $\V{y}$ are constructed by adding a noise vector to each ground-truth signal, such vector is sampled from $\calN(0,\sigma^2 \V{I}_{p})$, where $\sigma$ varies through experiments.
The ``true'' underlying dictionary $\V{D}_{TV}$ is shown in  \cref{fig:varying_sigma} (a), up to permuting its columns, and to rescaling with $\lambda$, which we compute in each experiment with a grid search solving \cref{eq:D_learning_opt_cost} \wrt $\lambda$.
Matrices $\hat{\V{D}}$ shown in this section have their columns sorted by magnitudes, to ease the comparison with $\V{D}_{TV}$.
We consider \emph{stochastic gradient descent} updates with batch size 64. Training set size: 640000, validation set size: 256. We adopt random white noise initialization with varying variance.

\begin{figure}
	\centering
	\includegraphics[width=.22\textwidth]{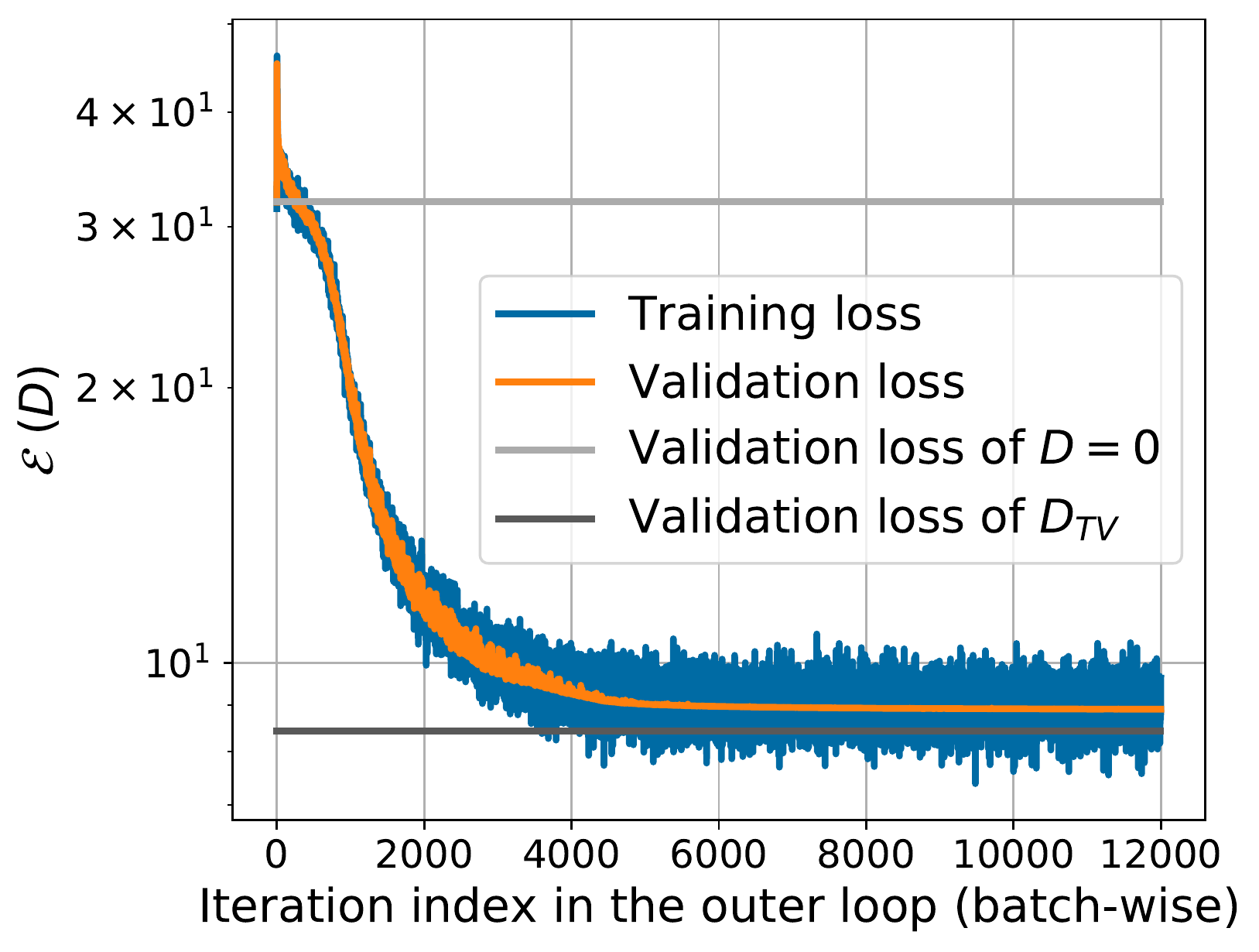}
	\hfill
	\includegraphics[width=.22\textwidth]{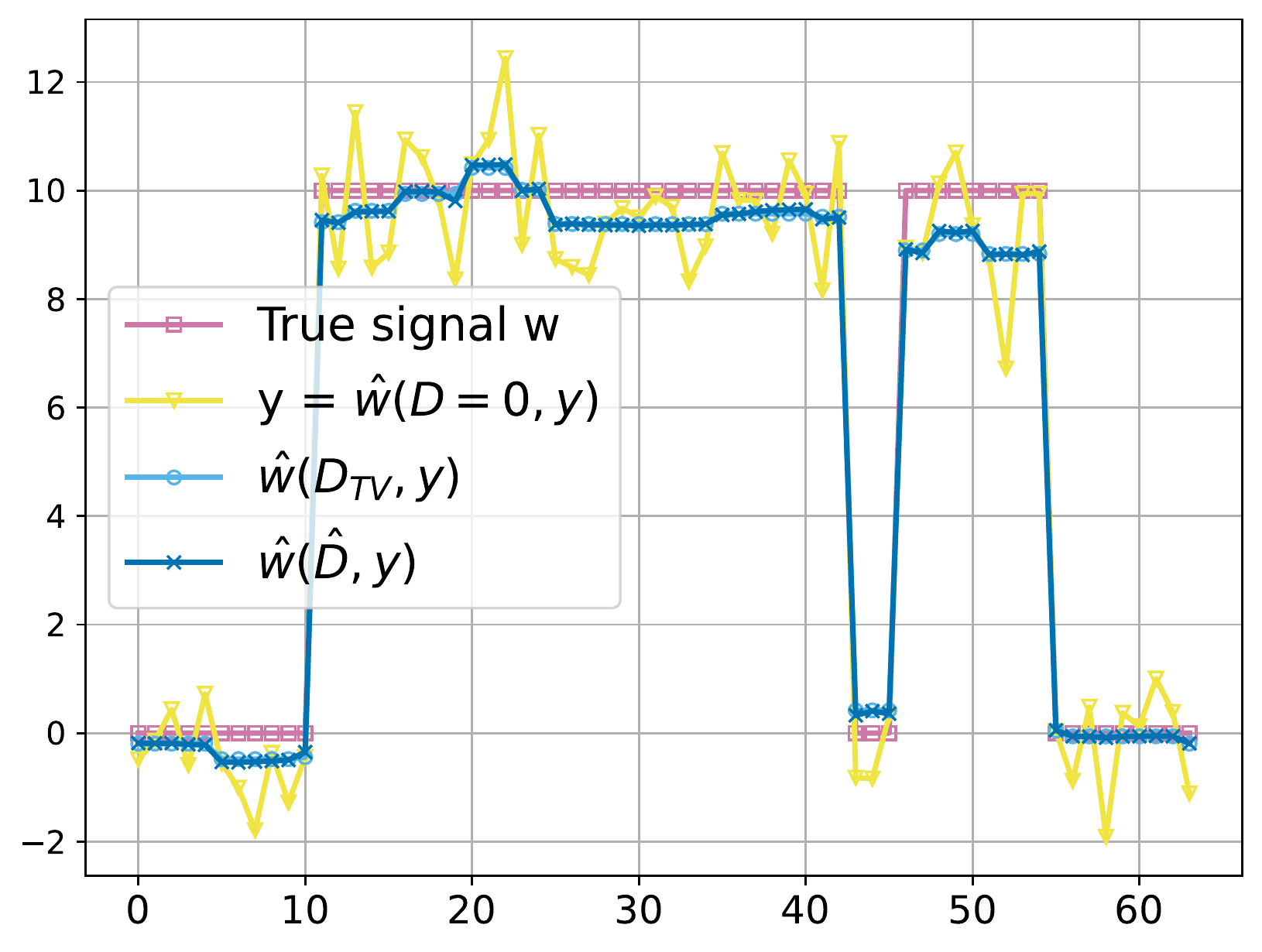}
	\caption{Performance curves of the learnt $\hat{\V{D}}$ in the case $\sigma=1$ in \cref{fig:varying_sigma}.
		Left: training and validation losses as in \cref{eq:D_learning_opt_cost} as a function of the iteration index. We plot the loss incurred by $\V{D}_{TV}$ and $\V{D}=\V{0}$ on the validation set.
		Right: Recovered signals with $\V{D}_{TV}$, $\V{D}=\V{0}$ and $\hat{\V{D}}$ as in \cref{eq:denoising_opt_cost} from a random signal $\V{y}$ from the validation set.
	}
	\label {Fig: train_val_loss}
\end{figure}

\textbf{Denoiser setup:} The value of $\eta_1$ in \cref{eq:dual_FB_update} is assigned automatically while guaranteeing convergence of $\V{z}_i$. This value is $\eta_1 = 0.95/\|\V{D}^\top \V{D}\|_2$, where $\|\V{D}^\top \V{D}\|_2$ is the Lipschitz constant of $\nabla \frac{1}{2} \big\Vert \V{Dz}- \V{y} \big\Vert_2^2$.
In all experiments, a tolerance threshold is used to determine the number of  iterations, specifically, we assume convergence if $\|\V{z}_{i+1}-\V{z}_i\|_\infty/\|\V{z}_i\|_\infty<10^{-4}$.

\subsection{Projection proposed for this case study}

We reduce the admissible set $\calC$ in our proposed algorithm \ref{alg} to dictionaries whose columns sum up to zero, \ie $\calC=\{\V{D} |\1_{p}^\top \V{D}  = \V{0}_m\}$.
This property is seen in the prior operator $\V{D}_{TV}$ our algorithm is expected to learn, and in the more general family of problems known as graph total variation \cite{berger2017graph}. This very simple prior greatly improves the results by filtering out many local minima. Different priors for other cases will be the goal of future investigations.
Referring to the $c$-th column of $\V{D}$ by $\V{D}[:,c]$, and by $\mean( 	{\V{D}}[:,c] )$ to the mean value of this column, we project $\V{D}_{t+1}$ as follows:
\begin{equation}
	\label{eq:projection}
	\V{D}_{t+1}[:,c] =
	{\V{D}}_{t+1}[:,c] - \mean\big({\V{D}}_{t+1}[:,c] \big),~ \forall c \in [m].
\end{equation}

\begin{figure}
	\centering
	\includegraphics[width=.45\textwidth]{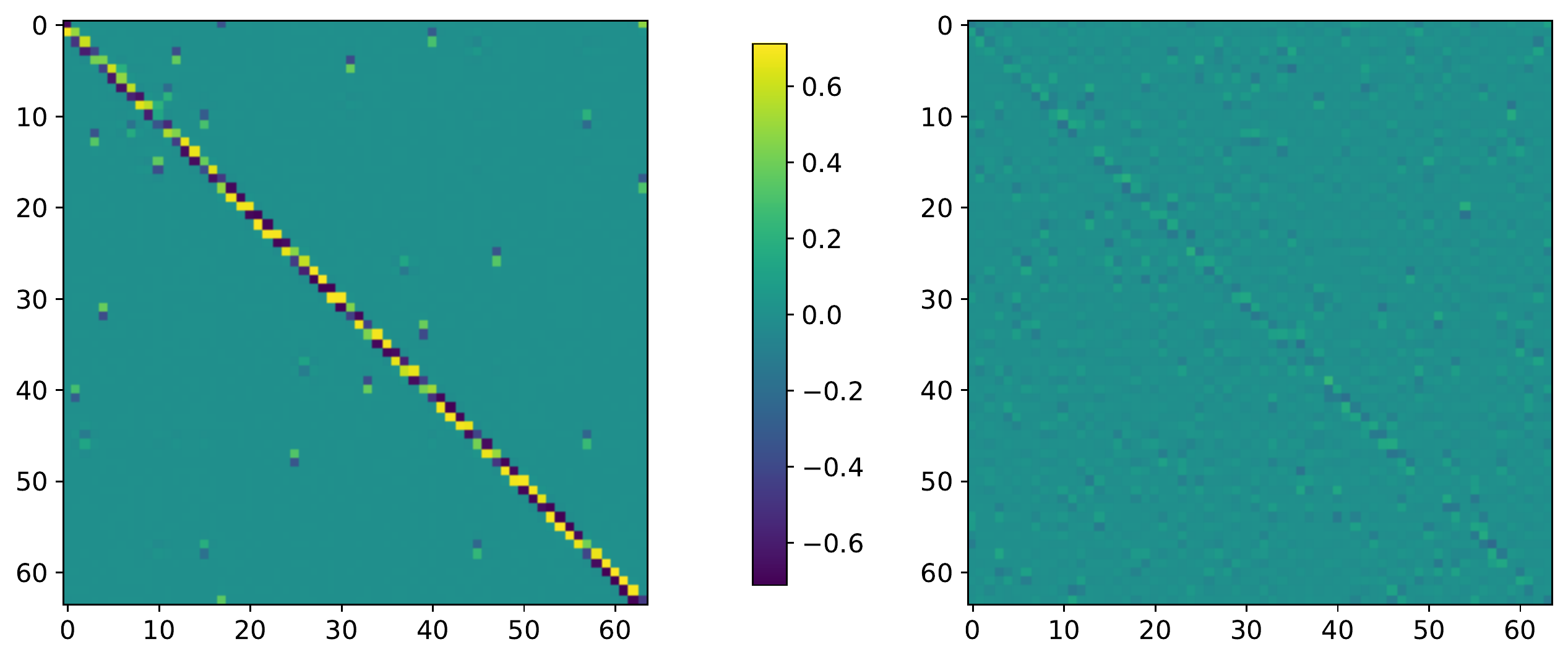}
	\caption{Showing the difference made by our proposed admissible set $\calC$. We plot $\hat{\V{D}}$ produced by our algorithm.  Left: with our proposed $\calC$.
		Right: with $\calC = \bbR^{p\times p}$.
		Dataset: Jumps in the dataset can randomly occur between any two values in $[0,10]$, thus they have random amplitudes. $\sigma=0.5$.}
	\label{fig:is_0colsum_useful}
\end{figure}

\begin{figure}
	\centering
	\includegraphics[width=.45\textwidth]{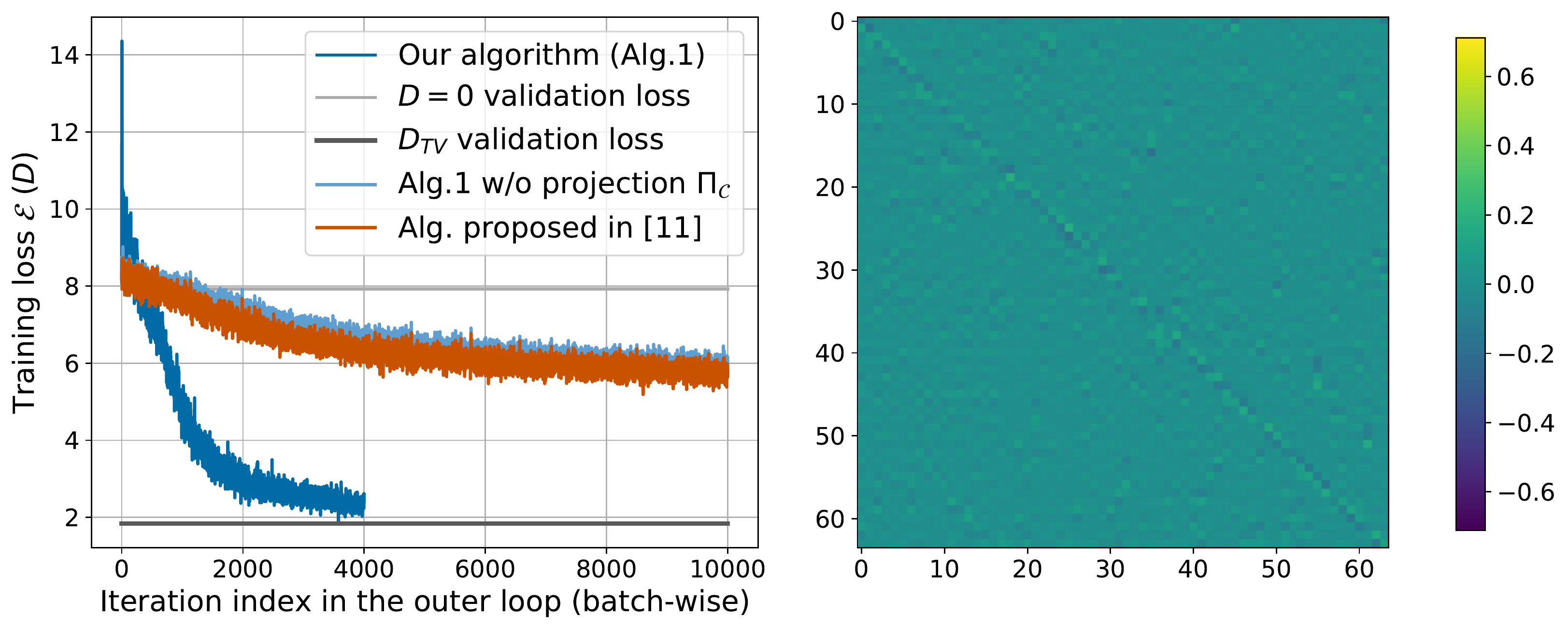}
	\caption{Benchmark against the unconstrained optimization with $\ell_1$-smoothening proposed in \cite{peyre2011learning}, with identical problem setting to \cref{fig:is_0colsum_useful}.
		Left: training loss curves of different algorithms.
		Right: $\hat{\V{D}}$ produced by the algorithm in \cite{peyre2011learning}.
		Dataset: same as in \cref{fig:is_0colsum_useful}.}
	\label{fig: Peyre_benchmark}
\end{figure}

In \cref{fig:is_0colsum_useful} and \cref{fig: Peyre_benchmark}(left), we show the essential role of the centering projection in our proposed algorithm. We run our algorithm twice: \emph{i)} with the projection; \emph{ii)} ignoring the projection step; on the same dataset, and we plot the learned dictionary $\hat{\V{D}}$ in both cases. Discontinuities in the used dataset are of random magnitudes (between any two values in [0,10]), and the noise has a standard deviation $\sigma = 0.5$. Our algorithm coupled with the projection successfully captures $\V{D}_{TV}$-like structure from the dataset, unlike when the projection is not considered, which shows its capability in this problem setting.

\subsection{Sensitivity \wrt to the noise level}

When our observations are noise-free, \ie $\sigma =0$ and $\V{y}=\V{w}$, it is clear $\V{D}=\V{0}$ is the optimal dictionary, as no regularization is required to retrieve the true signal.
Therefore, when $\sigma$ is small, we don't expect our algorithm to learn the structure in $\V{D}_{TV}$, since in such case, its magnitude is of the same level as numerical errors, see \cref{fig:varying_sigma} (b).

On the other hand, when the noise level is  high, the piecewise constant prior is degraded and is poorly seen in observations. As a result, the learnt dictionary is a distorted version of $\V{D}_{TV}$, as in \cref{fig:varying_sigma} (f).

In between, it is important to verify that our algorithm is stable, and can extract $\V{D}_{TV}$ out from data. As shown in \cref{fig:varying_sigma} (c-e), this is indeed true when $\sigma \in [0.75, 2]$, which spans a non-trivial range in SNR scale $[12,89]$, given that our signals $\V{w}$ take values in $\left\{ 0,10 \right\}$.

For $\sigma =1$, we show in \cref{Fig: train_val_loss} the evolution of training/validation losses as a function of the batch index through training, while benchmarking it against the same loss incurred when $\V{D}= \V{D}_{TV}$ with optimal $\lambda$ and $\V{D}= \V{0}$ as a reference. We also present the same benchmark but \wrt the denoising performance, as in \cref{eq:denoising_opt_cost}, applied on a signal from the validation set.

\subsection{Benchmark against the algorithm proposed in \cite{peyre2011learning} }

In \cref{fig: Peyre_benchmark}, we compare our output to the one produced by the unconstrained learning algorithm based on smoothening $\ell_1$ in \cite{peyre2011learning}. We fix the $\ell_1$-smoothening parameter $\epsilon=10^{-3}$. The dataset has discontinuities of random magnitudes in [0,10], and the noise standard deviation $\sigma = 0.5$, \ie same setup as in \cref{fig:is_0colsum_useful}. Unlike our projected gradient descent, the opponent algorithm fails to inspect the  $\V{D}_{TV}$ structure  from data.

\section{Conclusion and future work}
To alleviate the difficulty of deriving gradients to learn analysis-sparsity dictionaries, that are optimized through a bilevel problem, we proposed to make use of automatic differentiation, a technique shown to have high proficiency in machine learning and deep learning.
In the absence of a theoretical link between the automatic differentiation and the analytical one in our setting, our experiments on the piecewise constant signals reconstruction problem showed a proof-of-concept, and a promising methodology that can be applied in other problems setting.
For the same case study, we also incorporated a simple column-wise centering projection, which significantly increased the stability of the algorithm and the quality of the learned dictionary.

Further investigations can be done to generalize our algorithm to 2D signals, particularly the 2D piecewise constant signals reconstruction, while trying to combine other constraints. Moreover, it would be interesting to study the problem again when the dimensions of the dictionary are not given, for instance to learn graph incidence matrices.

\bibliography{refs}

\begin{thebibliography}{17}
\providecommand{\natexlab}[1]{#1}
\providecommand{\url}[1]{\texttt{#1}}
\expandafter\ifx\csname urlstyle\endcsname\relax
  \providecommand{\doi}[1]{doi: #1}\else
  \providecommand{\doi}{doi: \begingroup \urlstyle{rm}\Url}\fi

\bibitem[Baydin et~al.(2018)Baydin, Pearlmutter, Radul, and
  Siskind]{baydin2018automatic}
A.~G. Baydin, B.~A. Pearlmutter, A.~A. Radul, and J.~M. Siskind.
\newblock Automatic differentiation in machine learning: a survey.
\newblock \emph{J. Mach. Learn. Res.}, 18, 2018.

\bibitem[Beck and Teboulle(2009)]{beck2009fast}
A.~Beck and M.~Teboulle.
\newblock A fast iterative shrinkage-thresholding algorithm for linear inverse
  problems.
\newblock \emph{SIAM J. Imaging Sci.}, 2\penalty0 (1):\penalty0 183--202, 2009.

\bibitem[Berger et~al.(2017)Berger, Hannak, and Matz]{berger2017graph}
P.~Berger, G.~Hannak, and G.~Matz.
\newblock Graph signal recovery via primal-dual algorithms for total variation
  minimization.
\newblock \emph{IEEE Journal of Selected Topics in Signal Processing},
  11\penalty0 (6):\penalty0 842--855, 2017.

\bibitem[Chambolle and Pock(2020)]{chambolle2020learning}
A.~Chambolle and T.~Pock.
\newblock Learning consistent discretizations of the total variation.
\newblock 2020.

\bibitem[Chambolle et~al.(2010)Chambolle, Caselles, Cremers, Novaga, and
  Pock]{chambolle2010introduction}
A.~Chambolle, V.~Caselles, D.~Cremers, M.~Novaga, and T.~Pock.
\newblock An introduction to total variation for image analysis.
\newblock \emph{Theoretical foundations and numerical methods for sparse
  recovery}, 9\penalty0 (263-340):\penalty0 227, 2010.

\bibitem[Elad et~al.(2007)Elad, Milanfar, and Rubinstein]{elad2007analysis}
M.~Elad, P.~Milanfar, and R.~Rubinstein.
\newblock Analysis versus synthesis in signal priors.
\newblock \emph{Inverse problems}, 23\penalty0 (3):\penalty0 947, 2007.

\bibitem[Grabmeier and Kaltofen(2003)]{grabmeier2003computer}
J.~Grabmeier and E.~Kaltofen.
\newblock \emph{Computer Algebra Handbook: Foundations, Applications, Systems}.
\newblock Springer Science \& Business Media, 2003.

\bibitem[Jerrell(1997)]{jerrell1997automatic}
M.~E. Jerrell.
\newblock Automatic differentiation and interval arithmetic for estimation of
  disequilibrium models.
\newblock \emph{Computational Economics}, 10\penalty0 (3):\penalty0 295--316,
  1997.

\bibitem[{Mancera} and {Portilla}(2006)]{Mancera2006L0}
L.~{Mancera} and J.~{Portilla}.
\newblock L0-norm-based sparse representation through alternate projections.
\newblock In \emph{2006 International Conference on Image Processing}, pages
  2089--2092, 2006.

\bibitem[McCann and Ravishankar(2020)]{mccann2020supervised}
M.~T. McCann and S.~Ravishankar.
\newblock Supervised learning of sparsity-promoting regularizers for denoising.
\newblock \emph{arXiv preprint arXiv:2006.05521}, 2020.

\bibitem[Peyr{\'e} and Fadili(2011)]{peyre2011learning}
G.~Peyr{\'e} and J.~M. Fadili.
\newblock Learning analysis sparsity priors.
\newblock In \emph{Sampta'11}, pages 4--pp, 2011.

\bibitem[Ravishankar and Bresler(2015)]{ravishankar2015sparsifying}
S.~Ravishankar and Y.~Bresler.
\newblock Sparsifying transform learning with efficient optimal updates and
  convergence guarantees.
\newblock \emph{IEEE Transactions on Signal Processing}, 63\penalty0
  (9):\penalty0 2389--2404, 2015.

\bibitem[Rockafellar(1974)]{rockafellar1974conjugate}
R.~T. Rockafellar.
\newblock \emph{Conjugate duality and optimization}.
\newblock SIAM, 1974.

\bibitem[Sprechmann et~al.(2013)Sprechmann, Litman, Ben~Yakar, Bronstein, and
  Sapiro]{sprechmann2013supervised}
P.~Sprechmann, R.~Litman, T.~Ben~Yakar, A.~M. Bronstein, and G.~Sapiro.
\newblock Supervised sparse analysis and synthesis operators.
\newblock \emph{Advances in Neural Information Processing Systems},
  26:\penalty0 908--916, 2013.

\bibitem[Verma(2000)]{verma2000introduction}
A.~Verma.
\newblock An introduction to automatic differentiation.
\newblock \emph{Current Science}, pages 804--807, 2000.

\bibitem[Yaghoobi et~al.(2011)Yaghoobi, Nam, Gribonval, and
  Davies]{yaghoobi2011analysis}
M.~Yaghoobi, S.~Nam, R.~Gribonval, and M.~E. Davies.
\newblock Analysis operator learning for overcomplete cosparse representations.
\newblock In \emph{2011 19th European Signal Processing Conference}, pages
  1470--1474. IEEE, 2011.

\bibitem[Yaghoobi et~al.(2012)Yaghoobi, Nam, Gribonval, and
  Davies]{yaghoobi2012noise}
M.~Yaghoobi, S.~Nam, R.~Gribonval, and M.~E. Davies.
\newblock Noise aware analysis operator learning for approximately cosparse
  signals.
\newblock In \emph{2012 IEEE International Conference on Acoustics, Speech and
  Signal Processing (ICASSP)}, pages 5409--5412. IEEE, 2012.

\end{thebibliography}
\end{document}